\magnification=\magstep1
\hsize=16.5 true cm 
\vsize=23.6 true cm
\font\bff=cmbx10 scaled \magstep1
\font\bfff=cmbx10 scaled \magstep2
\font\bffg=cmbx10 scaled \magstep3

\parindent0cm
\def\cl{\centerline}           %
\def\bp{\bigskip}              %
\def\mp{\medskip}              %
\def\sp{\smallskip}            %
\def\bc{{\bf c}}               %
\font\boldmas=msbm10           %
\def\Bbb#1{\hbox{\boldmas #1}} %
\def\R{\Bbb R}                 %
\def\C{\Bbb C}                 %
\def\N{\Bbb N}                 %
\def\K{\Bbb K}                 %

\cl{\bffg Connected Hamel bases in Hilbert spaces}
\bp
\cl{\bfff Gerald Kuba}
\bp\bp

\bp\mp
{\bff 1. Introduction and main results}
\mp
The cardinal number (the {\it size}) of a set $\,S\,$
is denoted by $\,|S|\,$ and 
$\,\bc:=|\R|\,$. 
Throughout this note, $\,\K\in\{\R,\C\}\,$ and 
$\,X\,$ is an infinite-dimensional 
normed vector space over the field $\,\K\,$ 
and $\;\dim X\;$ denotes the transfinite
{\it algebraic dimension} of $\,X\,$ or, equivalently, 
the size of a basis of $\,X\,$. A {\it basis} is always an algebraic
basis (also called a Hamel basis)
and must not be confused with an orthonormal basis if $\,X\,$
is a Hilbert space. 
For example, $\;\dim X=\bc\;$ if $\,X\,$ is (real or complex) 
Hilbert space $\,\ell^2\,$, whereas every orthonormal basis
of $\,\ell^2\,$ is only a countable set.
An unorthodox reason why $\;\dim\ell^2=\bc\;$ is true is the following.
\mp
(1.1)\quad {\it There exists an injective, continuous function $\,f\,$ 
from $\,{]0,1[}\,$ into the unit sphere $\,||x||=1\,$ 
of Hilbert space $\,\ell^2\,$ such that the set $\,f(]0,1[)\,$ 
is linearly independent over $\,\K\,$.} 
\mp
One can easily settle (1.1) by defining 
$\;f(t)\,=\,(t^n\sqrt{1\!-\!t^2})_{n\in\N}\,$ for $\,0<t<1\,$.
(It is plain that $\,f\,$ is a continuous function. Actually, $\,f\,$ 
is $\,C^\infty\,$. By computing Vandermonde 
determinants as in [3], $\,f(]0,1[)\,$ is linearly independent.)
Of course, by applying Zorn's lemma, the 
linearly independent and pathwise connected set $\,f(]0,1[)\,$
of size $\,\bc\,$ can be expanded to a basis $\,B\,$ 
of $\,\ell^2\,$  such that $\,B\,$ lies in the unit sphere.
But there is no chance to accomplish that $\,B\,$ 
is {\it pathwise connected} as well. What we can accomplish is 
the following.
\mp
(1.2)\quad{\it There exists a basis of Hilbert space $\,\ell^2\,$
which is a connected and locally connected subset of the unit sphere.}
\mp
Our main goal is to prove a natural generalization of (1.2).
As usual, $\,o(X)\,$ denotes the total number of all open subsets of $\,X\,$. 
Naturally, $\;o(X)\geq |X|\geq\bc\;$ and $\;\dim X\leq|X|\,$.
Notice also that $\;\dim X=|X|\;$ if $\,|X|>\bc\,$.
(This is trivial in view of $\,|\C|\leq\bc\,$ and the fact 
that from an infinite set $\,B\,$ one can select only $\,|B|\,$
finite subsets.) Therefore, since the size of an infinite, {\it connected}
metric space cannot be smaller than $\,\bc\,$, we always have 
$\;\dim X=|X|\;$ if $\,X\,$ has a connected basis. 
It is possible that $\,o(X)>|X|\,$, see the remark below.
However, if $\,X\,$ is separable then $\,o(X)=|X|=\bc\,$. 
(There are precisely $\,\bc\,$ unions of sets taken from 
a countably infinite family of open balls.)
In particular, $\;o(X)=|X|=\dim X\;$ if $\,X\,$ is Hilbert space $\,\ell^2\,$.
Therefore, (1.2) is a consequence 
of the following theorem.
\mp
{\bf Theorem 1.} {\it If $\,o(X)=|X|=\dim X\,$ then there exists a
basis of $\,X\,$ which is a connected and locally connected 
subset of the unit sphere $\;||x||=1\;$ in $\,X\,$.}
\mp\sp
{\it Remark.} In Theorem 1
the case $\,o(X)=|X|=\dim X=\bc\,$ is the most important one
since it occurs if $\,X\,$ is one of the 
prominent separable Banach spaces $\,l^p\,$ and $\,L^p\,$
with $\;1\leq p <\infty\,$, including the Hilbert spaces $\,\ell^2\,$
and $\,L^2\,$. Nevertheless, 
in Theorem 1 the case $\;o(X)=|X|=\dim X>\bc\;$ is not vacuous
because if $\,X\,$ is a Hilbert space of weight $\,\kappa\,$
and $\,\kappa\,$ is any strong limit cardinal of countable cofinality 
then (see [2] (5.23)) $\,\kappa^{\aleph_0}=2^\kappa\,$ and hence 
(by [1] 4.1.H and 4.1.15) we have 
$\;o(X)=2^\kappa=\kappa^{\aleph_0}=|X|=\dim X\,$.
(However, if $\,\kappa^{\aleph_0}=\kappa\,$ for the weight 
$\,\kappa\,$ of $\,X\,$ 
then $\;o(X)=2^\kappa>\kappa=\kappa^{\aleph_0}=|X|=\dim X\,$.)
\vfill\eject
\bp
{\bff 2. Proof of Theorem 1}
\mp
Let $\,\varphi\,$ be the mapping $\;x\,\mapsto\,||x||^{-1}\cdot x\;$
from $\,X\setminus\{0\}\,$ onto the unit sphere 
$\;||x||=1\,$. Of course, $\,\varphi\,$ 
is continuous. Hence $\,\varphi(S)\,$ is connected whenever 
$\,S\,$ is a connected subset of $\,X\setminus\{0\}\,$.
Moreover, $\,\varphi\,$ is a retraction 
(i.e.~$\varphi\circ\varphi=\varphi\,$) 
and hence $\,\varphi(S)\,$ is locally connected whenever 
$\,S\,$ is a locally connected subset of $\,X\setminus\{0\}\,$
by applying [1] 2.4.E.(c) and 6.3.3.(d). 
It is evident that $\,\varphi\,$ restricted to a basis $\,B\,$
of $\,X\,$ is injective and $\,\varphi(B)\,$ is a basis of $\,X\,$
as well. Therefore Theorem 1 is a consequence of the following noteworthy 
theorem. 
\mp\sp
{\bf Theorem 2.} {\it If $\,o(X)=|X|=\dim X\,$ then there exists a
family $\,{\cal B}\,$ of 
bases of $\,X\,$ such that $\,|{\cal B}|=2^{|X|}\,$ and every 
element of $\,{\cal B}\,$  is a connected and locally connected and
dense subset of $\,X\,$.}
\mp
{\it Remark.} Obviously, the connected basis provided by 
Theorem 1 is {\it nowhere dense.}
\mp\sp
In order to prove Theorem 2 we need a crucial lemma about 
paths in convex sets.
In the following, if $\,A\subset X\,$ then $\,[A]\,$ is the 
linear subspace of $\,X\,$ generated by the 
vectors in the set $\,A\,$. Furthermore, if $\,L\,$ is a linear subspace 
of $\,X\,$ then $\;\hbox{\rm codim}\,L\;$ denotes the (possibly finite)
codimension of $\,L\,$ in $\,X\,$. (We have $\;\hbox{\rm codim}\,L=\kappa\;$
if and only if some basis $\,B\,$ of $\,X\,$ has a subset 
$\,B'\,$ with $\,[B']=L\,$ and $\;|B\setminus B'|=\kappa\,$.)
\mp                                                  
{\bf Lemma 1.} {\it If $\,U\not=\emptyset\,$ is a convex, open 
subset of $\,X\,$ and if $\,L\,$ is a linear subspace of $\,X\,$
with $\;\hbox{\rm codim}\, L\geq 2\;$
then $\;U\setminus L\;$ is dense in $\,U\,$
and pathwise connected.}
\mp
For the proof of Lemma 1 we need the following basic facts. 
\mp
(2.1) \quad {\it If the interior of a set $\,S\subset X\,$ is nonempty
then $\,[S]=X\,$.}
\mp
(2.2) \quad {\it If $\,U\,$  is a nonempty and open subset of $\,X\,$
and $\,L\,$ is a linear subspace of $\,X\,$ then 

\quad\qquad$\,\,U\setminus L\,$ contains a linearly independent set of size 
$\,\hbox{\rm codim}\,L\,.$}
\mp
Notice that (2.2) covers the case $\,L=X\,$ since
$\;\hbox{\rm codim}\,X\,=\,0\;$
and $\,\emptyset\,$ is linearly independent.
It is evident that (2.2) follows from (2.1), while
a proof of (2.1) is a piece of cake. (If 
$\,V\,$ is a nonempty open ball with center $\,v\,$
then $\,W\,=\,-v+V\,$ is an open ball about $\,0\,$ and hence 
$\;[V]=[W]\,\supset\,\K\cdot W\,=\,X\,$.)
\mp
{\it Proof of Lemma 1.} 
Let $\,L'\,$ be any algebraic complement of $\,L\,$
in $\,X\,$. Then for every $\,x\in X\,$ there exists a unique pair 
$\;(u,v)\in L'\times L\;$ such that $\;x=u+v\;$. 
Let $\,x_1\,$ and $\,x_2\,$ be distinct points in $\,U\setminus L\,$
and write $\;x_j=u_j+v_j\;$ with $\;u_j\in L'\;$ and $\;v_j\in L\;$
for $\,j\in\{1,2\}\,$. Since $\;x_1,x_2\not\in L\;$ we have 
$\;u_1\not=0\not=u_2\,$.
We distinguish two cases.
Firstly assume that the vectors
$\;u_1,u_2\;$ are linearly independent.
Then we can be sure that $\;tu_1+(1-t)u_2\not=0\;$ 
for every $\,t\in\R\,$. Therefore, if 
$\;\ell(x,y)\;$ denotes the straight line segment 
which connects $\,x\in X\,$ with $\,y\in X\setminus\{x\}\,$
then we can be sure that $\,\ell(x_1,x_2)\,$ is disjoint from $\,L\,$.
Since $\,U\,$ is convex, $\,\ell(x_1,x_2)\subset U\,$
and hence $\,\ell(x_1,x_2)\,$ is a path 
in the topological space 
$\;U\setminus L\;$ which connects the points $\,x_1\,$ and $\,x_2\,$. 
\sp
Secondly, assume that the vectors $\;u_1,u_2\;$ are not linearly independent.
Equivalently, $\;u_1=\lambda u_2\;$ for some (real or complex) scalar 
$\,\lambda\,$. Since the (possibly finite) algebraic dimension
of the vector space $\,L'\,$ is greater than 1, in view of (2.2) we can 
choose a point $\,y\,$ in $\;U\setminus[L\cup\{u_1\}]\,$.
We claim that $\;\ell(y,x_j)\cap L=\emptyset\;$ 
for $\,j\in\{1,2\}\,$. Then, since $\,U\,$ is convex, 
$\;\ell(x_1,y)\cup\ell(x_2,y)\;$ is a path 
in the topological space 
$\;U\setminus L\;$ which connects the points $\,x_1\,$ and $\,x_2\,$. 
\sp
Assume on the contrary that $\;\ell(y,x_j)\cap L\not=\emptyset\;$ 
for some $\,j\in\{1,2\}\,$. Then for some $\,t\in[0,1]\,$
the point $\;ty+(1-t)x_j\;$ lies in $\,L\,$. 
Write $\;y\,=\,u_3+v_3\;$ with $\,u_3\in L'\,$ and $\,v_3\in L\,$.
Then $\;(tu_3+(1-t)u_j)+(tv_3+(1-t)v_j)\,\in\,L\;$
and this is only possible if $\;tu_3+(1-t)u_j=0\,$.
From $\;tu_3+(1-t)u_j=0\;$ (and $\,u_j\not=0\,$) we derive $\,t\not=0\,$ 
and hence $\,u_3\,$ lies in $\,[\{u_j\}]=[\{u_1\}]\;$
and hence $\,u_3\,$ lies in $\;[L\cup\{u_1\}]\,$. 
Since trivially $\;v_3\in[L\cup\{u_1\}]\,$, we conclude 
that $\;y=u_3+v_3\;$ lies in $\;[L\cup\{u_1\}]\;$
and this is the desired contradiction. 
\sp
Thus in both cases 
the points $\,x_1\,$ and $\,x_2\,$ can be connected by a path 
in $\;U\setminus L\;$ and hence 
$\;U\setminus L\;$ is pathwise connected.
And (2.1) implies that $\,L\,$ cannot contain a nonempty open set and hence 
$\,U\setminus L\,$ must be dense in $\,U\,$, q.e.d.
\bp\mp
{\bff 3. Proof of Theorem 2}
\mp
Let the transfinite cardinal number $\,\kappa\,$ be the dimension of $\,X\,$
and regard $\,\kappa\,$ as an initial ordinal number.
So for every ordinal $\,\alpha<\kappa\,$ we have 
$\;|\{\,\beta\;|\;\beta<\alpha\,\}|<\kappa\,$.
Let $\,{\cal G}\,$ denote the family of all subsets $\,G\,$
of $\,X\,$ such that $\,\dim[G]=\kappa\,$
and $\;G\,=\,U_1\setminus U_2\;$ for some pair $\,U_1,U_2\,$
of open subsets of $\,X\,$. Trivially, $\;|{\cal G}|\leq o(X)=\kappa\,$. 
In view of (2.1) every nonempty open subset of $\,X\,$ and every nondegenerate
closed ball in $\,X\,$ lies in the family $\,{\cal G}\,$.
Consequently, $\;|{\cal G}|=o(X)=\kappa\,$. 
\mp
Since $\;|\{\,\alpha\;|\;\alpha<\kappa\,\}|=\kappa=|{\cal G}|\,$,
we can write $\;{\cal G}\;=\;\{\,G_\alpha\;|\;\alpha<\kappa\,\}\;$ 
(and we do not care whether the mapping $\,\alpha\mapsto G_\alpha\,$
is injective or not).
Let $\,\rho\,$ be a choice function defined 
on the nonempty subsets of $\,X\,$.
Thus $\,\rho(S)\in S\,$ 
for  every nonempty set $\,S\subset X\,$. 
\mp
Now, by induction we define points $\,x_\alpha\,$ in $\,X\,$ 
for all $\,\alpha<\kappa\,$.
Assume that for $\,\xi<\kappa\,$ a point 
$\,x_\alpha\,$ is already defined 
for every $\,\alpha<\xi\,$. (This assumption is vacuous if
$\,\xi=0\,$.)
Then define the point $\,x_\xi\,$ by
\medskip  
\centerline{$x_\xi\;:=\;\rho\big(\,G_\xi\setminus
[\{\,x_\alpha \;|\;\alpha<\xi\,\}]\,\big)\;$.}
\medskip
This definition is correct because the set 
$\;G_\xi\setminus [\{\,x_\alpha \;|\;\alpha<\xi\,\}]\;$
is nonempty. Indeed, $\;\dim[G_\xi]=\kappa\;$ 
implies that $\,G_\xi\,$ 
contains a linearly independent set $\,S\,$ of size $\,\kappa\,$,
whereas from $\;\dim[\{\,x_\alpha \;|\;\alpha<\xi\,\}]\,\leq\,
|\{\,\alpha \;|\;\alpha<\xi\,\}|<\kappa\;$
we conclude that this set $\,S\,$ cannot be a subset 
of $\;[\{\,x_\alpha \;|\;\alpha<\xi\,\}]\,$.
\medskip
In this way we obtain points $\;x_\alpha\;(\alpha<\kappa)\;$
where $\;x_\alpha\not=x_\beta\;$ whenever $\,\alpha<\beta<\kappa\,$.
In particular $\;A\,:=\,\{\,x_\alpha\;|\;\alpha<\kappa\,\}\;$
is a subset of $\,X\,$ with $\,|A|=\kappa\,$.
By construction, for every $\,\xi<\kappa\,$ the vector 
$\,x_\xi\,$ is linearly independent from all vectors 
$\;x_\alpha\;(\alpha<\xi)\,$. Therefore the whole set 
$\,A\,$ is linearly independent. 
\mp
Let $\,\Omega\,$ denote the set of all ordinals $\,\alpha<\kappa\,$ and let 
$\,\Omega_0\,$ be the set of all $\,\alpha\in\Omega\,$ such that 
$\,G_\alpha\,$ is 
an open ball of $\,X\,$.
Of course, $\,|\Omega_0|=\kappa\,$.
Fix one vector $\,x_\beta\in A\,$ with $\,\beta\not\in \Omega_0\,$
and choose for every $\,\alpha\in \Omega_0\,$ a scalar 
$\,\lambda_\alpha\not=0\,$
such that both points $\;x_\alpha\,$ and 
$\;\tilde x_\alpha\,:=\,x_\alpha+\lambda_\alpha x_\beta\;$ 
lie in the ball $\,G_\alpha\,$.
Let $\,{\cal F}\,$ be the family of all functions 
from $\,\Omega\,$ into $\,X\,$ where 
$\;f(\alpha)\in\{x_\alpha,\tilde x_\alpha\}\;$ 
for every $\,\alpha\in \Omega_0\,$
and $\;f(\alpha)=x_\alpha\;$ whenever $\,\alpha\not\in \Omega_0\,$.
Naturally, $\,|{\cal F}|=2^\kappa\,$. 
By construction,
$\,f(\Omega)\cap G\not=\emptyset\,$ for every $\,G\in{\cal G}\,$
and every $\,f\in{\cal F}\,$.
\mp
It is evident that  for 
every $\,f\in {\cal F}\,$ the linear subspace  
$\,[f(\Omega)]\,$ of $\,X\,$ coincides with $\,[A]\,$ and 
that $\,f(\Omega)\,$ is a basis of $\,[A]\,$.
By applying Zorn's lemma we can fix a subset $\,Y\,$ of $\,X\setminus[A]\,$
such that $\,Y\,$ is linearly independent and $\,[A\cup Y]=X\,$.
(If already $\,[A]=X\,$ then $\,Y=\emptyset\,$.)  
Consequently, $\;f(\Omega)\cap Y=\emptyset\;$
and $\;f(\Omega)\cup Y\;$ is a basis of the vector space $\,X\,$
for every $\,f\in{\cal F}\,$.
Trivially, $\,f(\Omega)\not=g(\Omega)\,$ whenever 
$\,f,g\in{\cal F}\,$ are distinct.
Therefore we can be sure that the size of 
\sp
\cl{$\;{\cal B}\,:=\,\{\,f(\Omega)\cup Y\;\,|\,\;f\in{\cal F}\,\}\;$}
\sp
equals $\,|{\cal F}|=2^\kappa\,$.
Now consider any basis $\,B\in{\cal B}\,$ of $\,X\,$.
Then, a fortiori, $\,B\,$ meets every set $\,G\in{\cal G}\,$
and in particular $\,B\,$ meets every nonempty open set, 
whence $\,B\,$ is a dense subset of $\,X\,$.
The proof is finished by verifying that $\,B\,$
is both connected and locally connected. 
Since $\,X\,$ and all open balls in $\,X\,$ are convex, 
it suffices to verify the following statement.
\mp
(3.1)\quad{\it If $\,U\subset X\,$ is convex and open then 
$\,B\cap U\,$ is connected.}
\mp
In order to verify (3.1) 
suppose on the contrary that $\,B\cap U\,$ is not connected.
Then there exist open sets $\,U_1,U_2\subset X\,$ 
such that $\;U_i\cap B\cap U\not=\emptyset\;$ for $\,i\in\{1,2\}\,$
and $\;B\cap U\,\subset\,U_1\cup U_2\;$ and 
$\;B\cap U\cap U_1\cap U_2=\emptyset\;$. 
Since $\,B\,$ is dense in $\,X\,$, from 
$\;B\cap U\cap U_1\cap U_2=\emptyset\;$
we infer $\;U\cap U_1\cap U_2=\emptyset\,$.
Therefore, since $\,U\,$ is connected, 
$\;G\;:=\;U\setminus(U_1\cup U_2)\;$ cannot be empty. 
\mp
Thus $\;U\setminus G\;$ is separated by the nonempty open sets 
$\,U\cap U_1\,$ and $\,U\cap U_2\,$ and hence $\,U\setminus G\,$ 
is not connected. Now, 
if $\,S\,$ is a dense and connected subset of a topological space
then every subset of the space containing $\,S\,$ is 
connected as well. Therefore, the subset 
$\,U\setminus[G]\,$ of the disconnected set $\,U\setminus G\,$
cannot be both connected and dense in $\,U\,$.    
Thus by Lemma~1 it is impossible 
that $\;\hbox{\rm codim}\,[G]\geq 2\,$.
Since $\,\dim X\,$ is infinite, from $\;\hbox{\rm codim}\,[G]< 2\;$ 
we conclude $\;\dim [G]=\dim X\;$ 
and hence $\,G\,$ lies in our family $\,{\cal G}\,$. 
Consequently, $\,B\cap G\not=\emptyset\,$. But this is impossible 
since $\;B\cap U\,\subset\,U_1\cup U_2\;$ and $\;G\subset U\,$.
Thus (3.1) is true and this concludes the proof of Theorem 2.
\bp
{\it Remark.} For $\;\varphi(x)=||x||^{-1}\cdot x\;$ we obviously 
have $\,\varphi(B_1)\not=\varphi(B_2)\,$ whenever 
$\,B_1,B_2\in{\cal B}\,$ are distinct.
Consequently, {\it if $\,o(X)=|X|=\dim X\,$ then there exist $\,2^{|X|}\,$ 
bases of $\,X\,$ which are connected and locally connected 
subsets of the unit sphere of $\,X\,$.}

\bp
\bigskip\bigskip
{\bff References} 
\medskip\smallskip
[1] Engelking, R.: {\it General 
Topology, revised and completed edition.} Heldermann 1989. 
\mp
[2] Jech, T.: {\it Set Theory.} 3rd ed. Springer 2002. 
\mp
[3] Kuba, G.: {\it Transfinite dimensions.} arXiv:2010.01983 [math.GM] 2020.

\bigskip
\bp\bp\bp
{\sl Author's address:} Institute of Mathematics, 

University of Natural Resources and Life Sciences, Vienna, Austria. 
\smallskip
{\sl E-mail:} {\tt gerald.kuba@boku.ac.at}

\end